\documentclass[11pt, leqno]{article}
\usepackage{amsmath, amsfonts, amsthm}
\begin{document}


\title{A new approach to  the Tarry-Escott problem}
\date{}
\author{Ajai Choudhry}

\maketitle

\theoremstyle{plain}

\newtheorem{lem}{Lemma}
\newtheorem{thm}{Theorem}

\begin{abstract} In this paper we describe a new method  of obtaining ideal solutions of the well-known Tarry-Escott problem, that is, the problem of   finding two distinct sets
of integers $\{x_1,\,x_2,\,\ldots,\,x_{k+1}\}$ and $\{y_1,\,y_2,\,\ldots,\,y_{k+1}\}$ 
 such that
$ \sum_{i=1}^{k+1}x_i^r=\sum_{i=1}^{k+1}y_i^r,\;\;\;r=1,\,2,\,\dots,\,k$, 
where $k$ is a given positive integer. When $k > 3$, only a limited number of parametric/ numerical ideal solutions of the Tarry-Escott problem are known. In this paper, by applying the new method mentioned above,  we  find several new parametric ideal  solutions of the problem when $k \leq 7$. The ideal solutions obtained by this new approach are more general and very frequently, simpler than the ideal solutions obtained by the earlier methods. We also obtain  new parametric solutions of certain diophantine systems that are closely related to the Tarry-Escott problem. These solutions are also more general and simpler than the solutions of these diophantine systems published earlier.
\end{abstract}



\renewcommand{\theequation}{\arabic{section}.\arabic{equation}}
\section{Introduction}\label{intro}
\hspace{0.25in} The Tarry-Escott problem (henceforth written briefly  as TEP)  of degree $k$ consists of finding two distinct sets
of integers $\{x_1,\,x_2,\,\ldots,\,x_s\}$ and $\{y_1,\,y_2,\,\ldots,\,y_s\}$ such that
\begin{equation}
\sum_{i=1}^sx_i^r=\sum_{i=1}^sy_i^r,\;\;\;r=1,\,2,\,\dots,\,k,
\label{tepsk}
\end{equation}
where $k$ is a given positive integer.

According to a well-known theorem of  Frolov  \cite[p.\  614]{Dor2}, the relations \eqref{tepsk} imply  that for any arbitrary integers $M$ and $K$,
\begin{equation}
\sum_{i=1}^s(Mx_i+K)^r=\sum_{i=1}^s(My_i+K)^r,\;\;\;r=1,\,2,\,\dots,\,k.
\end{equation}
Thus, if $x_i=a_i,\,y_i=b_i, \;i=1,\,2,\,\ldots,\,s$ is a solution of the diophantine system \eqref{tepsk}, then  $x_i=Ma_i+K,\;y_i=Mb_i+K,\;i=1,\,2,\,\dots,\,s$, 
is also a solution of \eqref{tepsk}, and all such solutions will be considered equivalent.

It follows from Frolov's  theorem that for each solution of the diophantine system \eqref{tepsk}, there is an equivalent one such that $\sum_{i=1}^sx_i=0=\sum_{i=1}^sy_i$ and the greatest common divisor of  the integers $x_i,\,y_i,\;i=1,\,2,\,\dots,\,s,$ is 1. This is known as the reduced form of the solution. If $x_i,\,y_i,\;i=1,\,2,\,\ldots,\,s$ is the reduced form of some solution of \eqref{tepsk}, then it is easy to see that the only other possible reduced form of the solution is $-x_i,\,-y_i,\;i=1,\,2,\,\ldots,\,s$, and so the reduced form of a solution is essentially unique.

We recall that a solution of \eqref{tepsk} is said to be a symmetric solution if it  satisfies  the additional  conditions  (if necessary, on re-arrangement) $x_i=-y_i,\;i=1,\,2,\,\dots,\,s,$ when $s$ is odd,  or the conditions  
$x_i=-x_{s+1-i},\; y_i=-y_{s+1-i},\;\;i=1,\,2,\,\dots,\,s/2,$ when $s$ is even. 
 Solutions that are equivalent to symmetric solutions are also considered to be symmetric. Solutions that are not symmetric are called nonsymmetric.

It is well-known   that for a non-trivial
solution of (\ref{tepsk}) to exist, we must have $s\geq (k+1)$ \cite[p.\  616]{Dor2}. 
 Solutions of (\ref{tepsk}) with the minimum possible value of $s$, that is, with $s=k+1,$ are known as ideal solutions of the problem. The problem of finding such ideal solutions  has attracted considerable attention. The complete ideal  solution of the TEP is known only when $k=2$ or $3$ \cite[pp.\  52,\,55-58]{Di1}. A limited number of parametric ideal solutions of the TEP are known when $4 \leq k \leq 7$ (\cite{Bor}, \cite{Che}, \cite{Cho0}, \cite{Cho1},  \cite{Cho2},  \cite{Dor1}, \cite[pp.\  41-54]{Glo}),  and infinitely many numerical ideal solutions are known when $k=8,\,9$ (\cite{Let}, \cite{Smy}) and $k=11$ \cite{Cho3}.  

In this paper we describe  a new method of finding ideal solutions of the TEP. As a first step, we will find solutions of \eqref{tepsk} for large values of $s$ such that $x_i,\,y_i$ are the terms of a certain number of arithmetic progressions, and  we then use these solutions of the TEP to obtain ideal solutions. We describe the  method in greater detail in Section~\ref{prelim}.  This  new  method  could easily be applied to obtain  the complete ideal solution of the TEP of degrees 2 and 3 but since the complete solution is already known in these cases (see \cite[pp.\  52, 55-58]{Di1}),  we omit these solutions. We apply the new method in Sections~\ref{tepdeg4}, \ref{tepdeg5}, \ref{tepdeg6} and  \ref{tepdeg7} to  obtain several ideal solutions of the TEP of degrees $4,\,5,\,6$ and $7$ respectively.  All of these ideal solutions  will be presented in the reduced form. The ideal solutions obtained in this paper are more general and very frequently, simpler than the ideal solutions obtained by the earlier methods.

In Section~\ref{relsys} we use the new method  to derive parametric solutions of some    diophantine systems that are closely related to the TEP. These solutions are also more general and simpler than the known results concerning  these diophantine systems. As an example, we obtain a simple parametric solution, in terms of polynomials of degree 2,  of the following diophantine system,
\begin{equation}
\begin{aligned}
\sum_{i=1}^7x_i^7&=\sum_{i=1}^7y_i^r,\;\;\;r=1,\,2,\,\dots,\,5,\\
\prod_{i=1}^7x_i&=\prod_{i=1}^7y_i.
\end{aligned} \label{tep5eqprod}
\end{equation}
Till now only one parametric solution of the diophantine system \eqref{tep5eqprod}, in terms of polynomials of degree 14, had been obtained \cite{Cho4}.

\section{A general method of obtaining new solutions of the Tarry-Escott problem}\label{prelim}
 \setcounter{equation}{0}
\hspace{0.25in} Our new approach to the TEP is based on a well-known lemma \cite[p.\ 615]{Dor2} which is given  below without proof.

\begin{lem}\label{lemTarry}
If there exist integers $x_i,\,y_i,\,i=1,\,2,\,\ldots,\,s,$ such that the relations \eqref{tepsk} are satisfied, then 
\begin{equation}
\sum_{i=1}^s\{x_i^r+(y_i+h)^r\}=\sum_{i=1}^s\{(x_i+h)^r+y_i^r\},\;\;\;r=1,\,2,\,\dots,\,k,\,k+1, \label{tepskh}
\end{equation}
where $h$ is an arbitrary integer.
\end{lem}

We will now describe   the  general method  adopted in this paper.  We will first  find a  solution of the diophantine system \eqref{tepsk} taking the numbers $x_i$ and $y_i$ as the terms of  three or more  arithmetic progressions and directly solving the resulting diophantine equations. We will not impose any upper limit on the number of terms $s$ on either side of \eqref{tepsk}  and, in fact,  the integer $s$ could be much larger than $k$.  If $x_i=a_i,\,y_i=b_i,\;i=1,\,2,\,\ldots,\,s$, is any solution of the diophantine system \eqref{tepsk}, we will write this solution briefly as,
\begin{equation}
a_1,\,a_2,\,\ldots,\,a_s \stackrel{k}{=}  b_1,\,b_2,\,\ldots,\,b_s.
\label{notation}
\end{equation}

Let us assume that we obtain a  solution of \eqref{tepsk} given by,
\begin{equation}
\alpha_1,\,\alpha_2,\,\ldots,\,\alpha_m,\,\beta_1,\,\beta_2,\,\ldots,\,\beta_n \stackrel{k}{=}\gamma_1,\,\gamma_2,\,\ldots,\,\gamma_m,\,\delta_1,\,\delta_2,\,\ldots,\,\delta_n,
\label{solgen}
\end{equation}
with $\alpha_1,\,\alpha_2,\,\ldots,\,\alpha_m$ and $\beta_1,\,\beta_2,\,\ldots,\,\beta_n$ being  two arithmetic progressions consisting of $m$ and $n$ terms respectively where $m$ and $n$ could be arbitrarily large integers, while   $\gamma_1,\,\gamma_2,\,\ldots,\,\gamma_m$ and $\delta_1,\,\delta_2,\,\ldots,\,\delta_n$ are similarly  the terms of two arithmetic progressions consisting of $m$ and $n$ terms respectively.  

If the common difference of all the four arithmetic progressions is the same, say $d$, we apply Lemma~\ref{lemTarry} to the solution \eqref{solgen} taking $h=d$, and thus obtain the solution,
\[\begin{array}{l}
\alpha_1,\,\alpha_2,\,\ldots,\,\alpha_m,\,\beta_1,\,\beta_2,\,\ldots,\,\beta_n,\,\gamma_2,\,\gamma_3,\,\ldots,\,\gamma_m+d,\,
\delta_2,\,\delta_3,\,\ldots,\,\delta_n+d\\
\stackrel{k+1}{=}
\alpha_2,\,\alpha_3,\,\ldots,\,\alpha_m+d,\,\beta_2,\,\beta_3,\,\ldots,\,\beta_n+d,\,\gamma_1,\,\gamma_2,\,\ldots,\,\gamma_m,\,
\delta_1,\,\delta_2,\,\ldots,\,\delta_n,
\end{array}
\]
and on cancelling out the common terms on both sides, we obtain the solution,
\[
\alpha_1,\,\beta_1,\,\gamma_m+d,\,\delta_n+d \stackrel{k+1}{=} \alpha_m+d,\,\beta_n+d,\,\gamma_1,\,\delta_1.
\]

As a further example, if in the  solution \eqref{solgen}, the arithmetic progressions  $\alpha_1,\,\alpha_2,\,\ldots,\,\alpha_m$ and $\gamma_1,\,\gamma_2,\,\ldots,\,\gamma_m$ have common difference $d_1$ while the arithmetic progressions $\beta_1,\,\beta_2,\,\ldots,\,\beta_n$ and $\delta_1,\,\delta_2,\,\ldots,\,\delta_n$ have common difference $d_2$ with $d_1 \neq d_2$, we will apply Lemma~\ref{lemTarry} twice in succession, first taking $h=d_1$ when we get, after cancellation of common terms,  the  solution,
\[
\begin{array}{l}
\alpha_1,\,\beta_1,\,\beta_2,\,\ldots,\,\beta_n,\, \gamma_m+d_1,\,\delta_1+d_1,\,\delta_2+d_1,\,\ldots,\,\delta_n+d_1\\
\;\;\stackrel{k+1}{=}
\alpha_m+d_1,\,\beta_1+d_1,\,\beta_2+d_1,\,\ldots,\,\beta_n+d_1,\gamma_1,\,\delta_1,\,\delta_2,\,\ldots,\,\delta_n,
\end{array}
\]
and to this solution, we again apply Lemma~\ref{lemTarry}, this time taking $h=d_2$, and we thus get, after cancellation of common terms, the following solution:
\[
\begin{array}{l}
\alpha_1,\,\beta_1,\,\gamma_m+d_1,\, \delta_1+d_1,\,\alpha_m+d_1+d_2,\,\beta_n+d_1+d_2,\,
\gamma_1+d_2,\,\delta_n+d_2\\
\;\;\stackrel{k+2}{=}\alpha_1+d_2,\,\beta_n+d_2,\,\gamma_m+d_1+d_2,\,\delta_n+d_1+d_2,\,
\alpha_m+d_1,\,\beta_1+d_1,\,\gamma_1,\,\delta_1.
\end{array}
\]

It would be observed that in both of the illustrative examples given above, the number of terms on either side of the final solution is independent of the number of terms on either side of \eqref{solgen}. If the initial solution \eqref{solgen} is in terms of a certain number of arbitrary parameters, so is the final solution, and by choice of parameters, we can further reduce the number of terms on either side of the final solution. We will see  examples of this type in  Sections~\ref{tepdeg4} and \ref{tepdeg6}. 

As already mentioned, we will  find solutions of \eqref{tepsk} with values of $x_i,\,y_i$ being given by terms of certain arithmetic progressions. For facility of computation, we will invariably choose the arithmetic progressions to consist of  an even number of terms of the type,
\begin{equation}
\begin{array}{c}
a-(2n-1)d,\; a-(2n-3)d,\;\ldots,\;a-3d,\;a-d,\\
a+d,\;a+3d,\; \ldots,\; a+(2n-3)d,\;a+(2n-1)d. \label{APadn}
\end{array}
\end{equation}
We will refer to the $2n$ terms given by \eqref{APadn}, and having common difference $2d$,  as the terms of the arithmetic progression $[a,\,n,\,d]$. 

We will denote by $S_k(a,\,n,\,d)$ the  sum of $k^{\rm th}$ powers of the terms  given by \eqref{APadn}, that is,
\begin{equation}
S_k(a,\,n,\,d)=\sum_{j=1}^{n}\{(a-(2j-1)d)^k+(a+(2j-1)d)^k\}.
\end{equation}
The following  formulae are readily obtained using any standard symbolic algebra software such as MAPLE or Mathematica, and will be used frequently:
\begin{align}
S_1(a,\,n,\,d)&=2na,\label{s1}\\
S_2(a,\,n,\,d)&=2na^2+2n(4n^2-1)d^2/3,\label{s2}\\
S_3(a,\,n,\,d)&=2na^3+2n(4n^2-1)ad^2,\label{s3} \\
S_4(a,\,n,\,d)&=2na^4+4n(4n^2-1)a^2d^2 \nonumber\\
& \quad \quad +2n(4n^2-1)(12n^2-7)d^4/15. \label{s4} 
\end{align}

We will now obtain solutions of \eqref{tepsk} with $x_i,\,y_i$, being given by the terms of certain arithmetic progressions  $[a_j,\,m_j,\,d_j]$ and $[b_j,\,n_j,\,d_j]$. The initial choice of the number of arithmetic progressions on either side of \eqref{tepsk} as well as the number of terms and the common difference of each arithmetic progression is to be made suitably so that the resulting diophantine equations can be solved. After obtaining a solution of \eqref{tepsk}, we will    apply Lemma~\ref{lemTarry} either once or twice, as necessary,  to obtain the desired  solutions of the TEP. In several cases the solutions obtained have $m_j,\,n_j$ as arbitrary parameters. In our initial assumption, $m_j$ and $n_j$ are necessarily integers. The final parametric solutions of the TEP are, however, a finite number of polynomial identities of finite degree, and since these identities  are true for all the infinitely many integer values of $m_j$ and $n_j$, it follows that they are also true for all rational values of the parameters $m_j$ and $n_j$.

We note that since all the equations of the diophantine system \eqref{tepsk} are homogeneous,  if $x_i=a_i,\,y_i=b_i$ is any solution  of \eqref{tepsk}, and $\rho$ is  any nonzero rational number, then $x_i=\rho a_i,\,y_i=\rho b_i$ is also a  solution of \eqref{tepsk}. Thus any rational solution of \eqref{tepsk} yields a solution in integers on multiplying through by a suitable constant. Hence it suffices to obtain rational solutions of \eqref{tepsk}, and in the parametric solutions that we obtain, the arbitrary parameters may be assigned any rational values.  For the sake of brevity, we will omit the factor $\rho$ while writing the solutions of \eqref{tepsk}.

\section{Ideal solutions of the Tarry-Escott problem of degree 4}\label{tepdeg4}
 \setcounter{equation}{0}

\hspace{0.25in} We note that the complete ideal symmetric solution of degree 4 has been given by Choudhry \cite{Cho1}. We will therefore obtain only nonsymmetric  ideal solutions of the TEP of degree 4, that is, solutions of the diophantine system,
\begin{equation}
\sum_{i=1}^5x_i^r=\sum_{i=1}^5y_i^r,\;\;\;r=1,\,2,\,3,\,4,
\label{tep4s5}
\end{equation}
such that the simplifying conditions for symmetric solutions are not satisfied.

 Till now, apart from  numerical solutions \cite{Bor, Shu}, parametric nonsymmetric solutions of \eqref{tep4s5} in terms of polynomials of degrees 3 and 8 have been published (\cite{Cho1}, \cite{Cho2}). In this paper we obtain infinitely many parametric solutions of the diophantine system \eqref{tep4s5} in terms of polynomials of degree 2 as well as some   multi-parameter solutions. 

In both   of the following  subsections, we first obtain a parametric solution of the diophantine system,
\begin{equation}
\sum_{i=1}^6X_i^r=\sum_{i=1}^6Y_i^r,\;\;\;r=1,\,2,\,3,\,4,
\label{tep4s6}
\end{equation}
and  then choose the parameters so that we get a relation $X_i=Y_j$ for some suitable $i,\,j$, and thus, we may cancel out these terms, and we then get a solution of the diophantine system \eqref{tep4s5}.

\subsection{}

We will first  obtain a solution of the diophantine system \eqref{tep4s5} by  finding a solution of \eqref{tepsk} with $k=2$ taking the numbers $x_i$ as the terms of the two arithmetic progressions $[a,\,m_1,\,d_1]$ and $[0,\,m_2,\,d_2]$ where $d_1 \neq d_2$, and the numbers $y_i$ as the terms of the arithmetic progression $[b,\,n,\,d_1]$.  As the number of terms on either side of \eqref{tepsk} must be the same,  the integers $m_1,\,m_2$ and $n$ must satisfy the condition,
\begin{equation}
  m_1+m_2=n.\label{3tep4eq1}
\end{equation}

The two equations of the diophantine system \eqref{tepsk}  are given by,
\begin{align}
S_1(a,\,m_1,\,d_1)+S_1(0,\,m_2,\,d_2)&=S_1(b,\,n,\,d_1), \\
S_2(a,\,m_1,\,d_1)+S_2(0,\,m_2,\,d_2)&=S_2(b,\,n,\,d_1),
\end{align}
and, on using the the formulae \eqref{s1} and \eqref{s2}, these equations may be written as,
\begin{equation}
m_1a=nb,\label{3tep4eq2}
\end{equation}
and
\begin{multline}
2m_1a^2+(2/3)m_1(4m_1^2-1)d_1^2+(2/3)m_2(4m_2^2-1)d_2^2\\
=2nb^2+(2/3)n(4n^2-1)d_1^2.\label{3tep4eq3}
\end{multline}
We now have to solve the simultaneous equations \eqref{3tep4eq1},   \eqref{3tep4eq2} and \eqref{3tep4eq3}. 

The complete solution of \eqref{3tep4eq1} and \eqref{3tep4eq2} is given by,
\begin{equation}
a=(m_1+m_2)r,\quad b=m_1r,\quad n=m_1+m_2, \label{3tep4valabn}
\end{equation}
where $r$ is an arbitrary rational parameter, and on substituting these values of $a,\,b$ and $n$ in \eqref{3tep4eq3}, we get a quadratic equation in $d_1,\,d_2$ and $r$ whose complete solution in rational numbers is readily obtained and is given by,
\begin{equation}
\begin{aligned}
d_1& = -\rho\{(12m_1^2+12m_1m_2+4m_2^2-1)p^2-(8m_2^2-2)pq+(4m_2^2-1)q^2\},\\
 d_2& =\rho\{(12m_1^2+12m_1m_2+4m_2^2-1)p^2-(24m_1^2+24m_1m_2+8m_2^2-2)pq\\
& \quad \quad +(4m_2^2-1)q^2\},\\
 r& = \rho\{(24m_1^2+24m_1m_2+8m_2^2-2)p^2-(8m_2^2-2)q^2\},
\end{aligned}
\label{3tep4valdr}
\end{equation}
where $p$ and $q$ are arbitrary integer parameters while $\rho$ is an arbitrary nonzero rational parameter.

We now have a solution of \eqref{tepsk} with $k=2$ in which the numbers $x_i$ consist of two arithmetic progressions whose common differences are $2d_1$ and $2d_2$ while the numbers $y_i$ consist of the terms of a single  arithmetic progression whose common difference is $2d_1$. We  now apply Lemma~\ref{lemTarry} twice, in succession, taking $h=2d_1$ and $2d_2$ respectively, and obtain a solution of the diophantine system \eqref{tep4s6} which is given by,
\begin{equation}
\begin{aligned}
X_1&=a-(2m_1-1)d_1,\quad &X_2&= a+(2m_1+1)d_1+2d_2, \\
X_3&=-(2m_2-1)d_2,\quad &X_4&= (2m_2+1)d_2+2d_1,\\
X_5&= b+(2n+1)d_1,\quad &X_6&= b-(2n-1)d_1+2d_2,\\
Y_1&=a+(2m_1+1)d_1,\quad &Y_2&= a-(2m_1-1)d_1+2d_2,\\
Y_3&= (2m_2+1)d_2,\quad &Y_4&= -(2m_2-1)d_2+2d_1,\\
Y_5& =b-(2n-1)d_1, \quad &Y_6&=b+(2n+1)d_1+2d_2,
\end{aligned}
\label{sol3tep4XY}
\end{equation}
with the values of $a,\,b,\,d_1,\,d_2$ and $n$ being given by \eqref{3tep4valabn} and  \eqref{3tep4valdr} where $p,\,q,\,m_1$ and $m_2$ are arbitrary parameters.

 We can choose the parameters in the solution \eqref{sol3tep4XY} of \eqref{tep4s6}  in several ways so that one  pair of terms, one on each side of this solution of \eqref{tep4s6},  cancels out and we may thus obtain several solutions of \eqref{tep4s5}. As an example, on taking $p = 2m_2+1,\, q = 2m_1+2m_2+1,$ in the solution \eqref{sol3tep4XY}, we get $X_6=Y_4$, and we thus   obtain a solution of \eqref{tep4s5} which is given in the reduced form by

\begin{equation}
\begin{aligned}
x_1 & = 56m_1^2m_2+60m_1m_2^2+12m_2^3+44m_1^2+70m_1m_2+24m_2^2\\
& \quad \quad \quad +10m_1+9m_2,\\
 x_2 & = -24m_1^2m_2+12m_2^3-36m_1^2-40m_1m_2-6m_2^2-10m_1-6m_2,\\
 x_3 & = 56m_1^2m_2+80m_1m_2^2+32m_2^3-16m_1^2+20m_1m_2+24m_2^2+4m_2,\\ 
x_4 & = -64m_1^2m_2-80m_1m_2^2-28m_2^3-16m_1^2-60m_1m_2-36m_2^2\\
& \quad \quad \quad-10m_1-11m_2,\\
 x_5 & = -24m_1^2m_2-60m_1m_2^2-28m_2^3+24m_1^2+10m_1m_2-6m_2^2\\
& \quad \quad \quad+10m_1+4m_2,
\end{aligned}
\label{sol3tep4s5x}
\end{equation}
 and 
\begin{equation}
\begin{aligned}
y_1 & = -24m_1^2m_2+12m_2^3+24m_1^2+40m_1m_2+24m_2^2\\
&\quad \quad \quad +10m_1+9m_2,\\ 
y_2 & = 56m_1^2m_2+60m_1m_2^2+12m_2^3-16m_1^2-10m_1m_2-6m_2^2\\
& \quad \quad \quad-10m_1-6m_2,\\
 y_3 & = -64m_1^2m_2-80m_1m_2^2-28m_2^3-16m_1^2-20m_1m_2\\
&\quad \quad \quad-6m_2^2+4m_2,\\ 
y_4 & = 56m_1^2m_2+80m_1m_2^2+32m_2^3+44m_1^2+60m_1m_2+24m_2^2\\
& \quad \quad \quad +10m_1+4m_2,\\ 
y_5 & = -24m_1^2m_2-60m_1m_2^2-28m_2^3-36m_1^2-70m_1m_2-36m_2^2\\
& \quad \quad \quad -10m_1-11m_2,
\end{aligned}
\label{sol3tep4s5y}
\end{equation}
where $m_1$ and $m_2$ are arbitrary parameters. By successively assigning  distinct fixed numerical values to $m_2$ in the above solution, we obtain infinitely many  solutions of \eqref{tep4s5} in terms of polynomials of degree 2.  

As a numerical example, taking $m_1=1,\, m_2=1$ in \eqref{sol3tep4s5x} and \eqref{sol3tep4s5y} yields the following nonsymmetric ideal solution of the TEP of degree 4:
\begin{equation}
57,\,-22,\,40,\,-61,\,-14 \stackrel{4}{=} 19,\,16,\,-42,\,62,\,-55. \label{tep4ex1}
\end{equation}

\subsection{} 

To find a second  solution of the diophantine system \eqref{tep4s5}, we  first find a solution of \eqref{tepsk} with $k=3$ taking the numbers  $x_i,\,y_i$ as the terms of six arithmetic progressions with the same common difference $2d$. Specifically, we  take the numbers $x_i$ as  the terms of the arithmetic progressions $[a_j,\, n_j,\, d],\; j=1,\,2,\,3,$  and  the numbers  $y_i$ as  the terms of the arithmetic progressions $[b_j,\, n_j,\, d],\; j=1,\,2,\,3$.

We now have to solve the following three equations obtained by taking $r=1,\,2$ and $3$  in \eqref{tepsk}, and using the formulae \eqref{s1}, \eqref{s2} and \eqref{s3}:
\begin{align}
\sum_{i=1}^32n_ia_i&=\sum_{i=1}^32n_ib_i, \label{tep4eq1}\\
 \sum_{i=1}^32n_ia_i^2&=\sum_{i=1}^32n_ib_i^2, \label{tep4eq2}
\end{align}
and
\begin{multline}
\sum_{i=1}^32n_ia_i^3+\{\sum_{i=1}^32n_i(4n_i^2-1)a_i\}d^2\\
=\sum_{i=1}^32n_ib_i^3+\{\sum_{i=1}^32n_i(4n_i^2-1)b_i\}d^2. \label{tep4eq3}\\
\end{multline}

 We will  solve equations \eqref{tep4eq1}, \eqref{tep4eq2}  together with the following two equations,
\begin{align}
\sum_{i=1}^32n_ia_i^3&=\sum_{i=1}^32n_ib_i^3, \label{tep4eq3a}\\
\sum_{i=1}^32n_i(4n_i^2-1)a_i&=\sum_{i=1}^32n_i(4n_i^2-1)b_i, \label{tep4eq3b}
\end{align}
when \eqref{tep4eq3} will be identically satisfied for all values of $d$.

We  take $b_3=0$ for simplicity, and solve Eqs.~\eqref{tep4eq1} and \eqref{tep4eq3b} for $b_1,\,b_2$ to get,
 \begin{equation}
\begin{aligned}
b_1&=\{n_1(n_1^2-n_2^2)a_1-n_3(n_2^2-n_3^2)a_3\}/\{n_1(n_1^2-n_2^2)\},\\
b_2&=\{n_2(n_1^2-n_2^2)a_2+n_3(n_1^2-n_3^2)a_3\}/\{n_2(n_1^2-n_2^2)\}.
\end{aligned}
\label{tep4valb}
\end{equation}
Substituting these values of $b_i,\,i=1,\,2$, in \eqref{tep4eq2}, transposing all terms to one side and removing the factor $2n_3a_3$, we get,
\begin{multline}
2(n_1-n_2)(n_2^2-n_3^2)n_1n_2a_1-2(n_1-n_2)(n_1^2-n_3^2)n_1n_2a_2\\
+(n_1-n_3)(n_2-n_3)(n_1^3-n_1^2n_2+n_1^2n_3-n_1n_2^2\\
-3n_1n_2n_3-n_1n_3^2+n_2^3+n_2^2n_3-n_2n_3^2-n_3^3)a_3=0. \label{tep4eq2a}
\end{multline}

We now solve \eqref{tep4eq2a} to get,
\begin{multline}
a_1=\{2(n_1-n_2)(n_1^2-n_3^2)n_1n_2a_2-(n_1-n_3)(n_2-n_3)(n_1^3-n_1^2n_2\\
+n_1^2n_3-n_1n_2^2
-3n_1n_2n_3-n_1n_3^2+n_2^3+n_2^2n_3-n_2n_3^2-n_3^3)a_3\}\\
\times \{2n_1n_2(n_1-n_2)(n_2^2-n_3^2)\}^{-1}. \label{tep4vala1}
\end{multline}

Substituting the values of $b_1,\,b_2$ and $a_1$ given by \eqref{tep4valb} and \eqref{tep4vala1} in \eqref{tep4eq3a},  transposing all terms to one side and removing the factor $n_3(n_1^2-n_3^2)a_3$, we get the following quadratic equation in $a_2$ and $a_3$:
\begin{multline}
12n_1^2n_2^2(n_1^2-n_2^2)^2a_2^2-12n_1^2n_2(n_1^2-n_2^2)(n_2-n_3)\\
\times (n_1^2-n_2^2-n_2n_3-n_3^2)a_2a_3+(n_2-n_3)^2(3n_1^6-5n_1^4n_2^2\\
-4n_1^4n_2n_3-5n_1^4n_3^2+n_1^2n_2^4+2n_1^2n_2^3n_3+5n_1^2n_2^2n_3^2+2n_1^2n_2n_3^3\\
+n_1^2n_3^4+n_2^6+2n_2^5n_3-n_2^4n_3^2-4n_2^3n_3^3-n_2^2n_3^4+2n_2n_3^5+n_3^6)a_3^2=0. \label{tep4eq3c}
\end{multline}
Eq.~\eqref{tep4eq3c} will have a rational solution for $a_2$ and $a_3$ if, and only if, its discriminant given by,
\begin{multline}
48n_1^2n_2^2(n_1^2-n_2^2)^2(n_2^2-n_3^2)^2(n_1+n_2+n_3)(n_1+n_2-n_3)\\
\times (-n_1+n_2+n_3)(n_1-n_2+n_3),\label{tep4dis}
\end{multline}
becomes a perfect square. We must therefore choose $n_i,\;i=1,\,2,\,3$, such that the function 
\begin{equation}
\psi(n_1,\,n_2,\,n_3)=3(n_1+n_2+n_3)(n_1+n_2-n_3)(-n_1+n_2+n_3)(n_1-n_2+n_3),\label{tep4dis1}
\end{equation}
is a perfect square.

For $\psi(n_1,\,n_2,\,n_3)$ to be a perfect square, there must exist integers $f$ and $g$ such that  
\begin{equation}
3(n_1+n_2+n_3)(n_1+n_2-n_3)f^2=(-n_1+n_2+n_3)(n_1-n_2+n_3)g^2,
\end{equation}
and it further follows that there must exist integers $u$ and $v$ such that 
\begin{equation}
\begin{aligned}
3fv(n_1+n_2-n_3) = gu(n_1-n_2+n_3), \\
(n_1+n_2+n_3)fu = (-n_1+n_2+n_3)gv.
\end{aligned}
\label{condn}
\end{equation}

On solving the two linear equations \eqref{condn} for $n_1,\,n_2,\,n_3$,  we get,
\begin{equation}
\begin{aligned}
n_1 &= fgu^2+(3f^2-g^2)uv-3fgv^2,\\
 n_2 &= -(3f^2+g^2)uv,\\
 n_3 &= -fg(u^2+3v^2),
\end{aligned}
\label{tep4valn}
\end{equation}
where $f,\,g,\,u$ and $v$ are arbitrary parameters. We now substitute these values of $n_i$ in Eq.~\eqref{tep4eq3c}  and then solve it to  get the following solution for $a_2$ and $a_3$:

\begin{equation}
\begin{aligned}
a_2 &= (fu-gv)(-gu+3fv)\{(3f^2+4fg-3g^2)u^2+\\
& \quad (12f^2-24fg+4g^2)uv+(-27f^2+12fg+3g^2)v^2\},\\
 a_3& = 2(3f^2+g^2)(gu^2+6fuv-3gv^2)(fu^2-2guv-3fv^2).
\end{aligned}
\label{tep4vala23}
\end{equation}

Using \eqref{tep4vala23} and the values of $n_i$ already obtained, we get the value of $a_1$ from \eqref{tep4vala1}, and finally,  we obtain the values of $b_1$ and $b_2$ from \eqref{tep4valb}. These values of $a_1,\,b_1,\,b_2$ are given by,
\begin{equation}
\begin{aligned}
a_1&= u\{2fgu+(3f^2-g^2)v\}\{(3f^2+g^2)u^2+(12f^2-4g^2)uv\\
& \quad \quad -(27f^2+24fg+9g^2)v^2\},\\
b_1&= v\{(3f^2-g^2)u-6fgv\}\{(9f^2+8fg+3g^2)u^2\\
& \quad \quad \quad \quad +(12f^2-4g^2)uv-(9f^2+3g^2)v^2\},\\
b_2& = (fu+gv)(gu+3fv)\{(9f^2-4fg-g^2)u^2\\
& \quad \quad \quad \quad+(12f^2-24fg+4g^2)uv-(9f^2+12fg-9g^2)v^2\}.
\end{aligned}
\label{tep4vala1b12}
\end{equation}

We now have  a solution of \eqref{tepsk} with $k=3$ and with the numbers  $x_i, \,y_i$  consisting  of three arithmetic progressions whose common difference is $2d$. We can thus  apply Lemma~\ref{lemTarry} taking $h=2d$  to    obtain the following  solution of \eqref{tep4s6}:
\begin{equation}
\begin{aligned}
X_i&=a_i-(2n_i-1)d,\quad &X_{i+3}&= b_i+(2n_i+1)d,\;i=1,\,2,\,3,\\
Y_i&=a_i+(2n_i+1)d,\quad &Y_{i+3}&=b_i-(2n_i-1)d,\;i=1,\,2,\,3,
\end{aligned}
\label{2tep4valXY}
\end{equation}
where $d$ is an arbitrary parameter, the values of $n_1,\,n_2,\,n_3,\,a_1,\,a_2,\,a_3,\,b_1,\,b_2$ are given in terms of arbitrary parameters $f,\,g,\,u$ and $v$ by \eqref{tep4valn}, \eqref{tep4vala23}, \eqref{tep4vala1b12},  and $b_3=0.$ 

We note that, in the above  solution, $d$ occurs only in the first degree.  To obtain ideal solutions of the TEP of degree 4, we simply choose $d$ such that two of the terms, one on each side, become equal, and thus they can be cancelled out. We thus obtain a four-parameter ideal solution of the TEP of degree 4. There are 36 ways to  choose the pair of terms to be cancelled and we can thus obtain several distinct four-parameter solutions of our problem. 

The solutions obtained above are too cumbersome to be written down explicitly. Accordingly we give below  an example of  just one nonsymmetric solution obtained as described above by choosing  $d$ such that $X_1=Y_6$ and then taking $f=2,\,g=1$. This   solution of the diophantine system \eqref{tep4s5} is in  terms of two parameters $u$ and $v$ and  denoting the polynomial $c_0u^n+c_1u^{n-1}v+c_2u^{n-2}v^2+\cdots+c_nv^n$ by $[c_0,\,c_1,\,c_2,\,\ldots,\,c_n]$, this two-parameter solution may be  expressed in reduced form as follows: 
\begin{equation}
\begin{aligned}
x_1 &= [-98,\,-1007,\,1804,\,687,\,774], 
&x_2 &= [2,\,858,\,-2396,\,-1518,\,1314],\\
 x_3 &=[-128,\,308,\,1404,\,132,\,144],
 &x_4 &= [122,\,853,\,-1056,\,-753,\,-1206],\\
 x_5 &= [102,-1012,244,1452,-1026],
 &y_1 &= [102,1133,-2096,-693,1314],\\
 y_2 &= [122,-642,1804,1782,-1206],
 &y_3 &= [-128,-572,-56,2772,144],\\
 y_4 &= [2,-407,1404,-2013,-1026],
 &y_5 &= [-98,488,-1056,-1848,774].
\end{aligned}
\label{tep4valxy}
\end{equation}

As a numerical example, taking $u=2,\,v=-1$ in \eqref{tep4valxy} yields  the nonsymmetric solution,
\begin{equation}
2184,\,-2011,\,164,\,-1466,\,1129 \stackrel{4}{=} -2186,\,1589,\,-516,\,1984,\,-871. \label{tep4ex2}
\end{equation}

\section{Ideal solutions of the Tarry-Escott problem of degree 5}\label{tepdeg5}
 \setcounter{equation}{0}

\hspace{0.25in} Ideal  solutions of the TEP of degree 5 satisfy the system of equations,
\begin{equation}
\sum_{i=1}^6x_i^r=\sum_{i=1}^6y_i^r,\;\;\;r=1,\,2,\,\dots,\,5.
\label{tep5s6}
\end{equation}

Symmetric solutions of \eqref{tep5s6}  satisfy the additional conditions, 
\begin{equation}
x_4=-x_3,\,x_5=-x_2,\,x_6=-x_1,\,y_4=-y_3,\,y_5=-y_2,\,y_6=-y_1, \label{tep5condsym}
\end{equation}and under these  conditions, the diophantine system \eqref{tep5s6} reduces to the system of equations,
\begin{align}
x_1^2+x_2^2+x_3^2&=y_1^2+y_2^2+y_3^2, \label{tep5eq2}\\
x_1^4+x_2^4+x_3^4&=y_1^4+y_2^4+y_3^4.\label{tep5eq4}
\end{align}

The simultaneous diophantine equations \eqref{tep5eq2} and \eqref{tep5eq4} are solved very easily  if we impose   the additional condition $\sum_{i=1}^3x_i=\sum_{i=1}^3y_i$, since on taking $x_3=-x_1-x_2,\;y_3=-y_1-y_2$, both equations \eqref{tep5eq2} and \eqref{tep5eq4} reduce to the single quadratic equation $x_1^2+x_1x_2+x_2^2=y_1^2+y_1y_2+y_2^2$ whose complete solution is readily obtained. It is more interesting to find solutions of Eqs.~\eqref{tep5eq2} and \eqref{tep5eq4} together with the conditions,
\begin{equation}
\pm x_1 \pm x_2 \pm x_3 \neq \pm y_1 \pm y_2 \pm y_3,\quad \pm x_1 \pm x_2 \pm x_3 \neq 0. \label{lincondneq}
\end{equation}
 Only a limited number of parametric solutions of this type, in terms of polynomials of degrees 2 and 4, have been published (\cite{Che}, \cite[p.\  711]{Di2}, \cite[p.\  49]{Glo}). 

In Section 4.1 we obtain parametric symmetric solutions of the diophantine system \eqref{tep5s6} that yield much more general parametric solutions of Eqs.~\eqref{tep5eq2} and \eqref{tep5eq4} satisfying the conditions \eqref{lincondneq}. In Section 4.2 we  obtain a parametric nonsymmetric solution of \eqref{tep5s6} and show how more such solutions can be obtained. 

\subsection{Symmetric solutions of the TEP of degree 5}

\hspace*{0.25in} We will first obtain numerical   ideal symmetric solutions by directly solving equations \eqref{tep5eq2} and \eqref{tep5eq4} together with an additional condition namely, the numbers $x_i, i=1,\,2,\,3$ are in arithmetic progression, and so also are  the  numbers $y_i, i=1,\,2,\,3$. 
 We will thereafter obtain two parametric ideal symmetric solutions  by the general method of first finding a solution of the diophantine system \eqref{tepsk} such that $x_i,\,y_i, \;i=1,\,2,\,\dots,\,s$, are the terms of arithmetic progressions.

\subsubsection{} To  solve equations \eqref{tep5eq2} and \eqref{tep5eq4} directly, we   write, 
\begin{equation}
x_1=a-d_1,\;\;x_2=a,\;\;x_3=a+d_1,\;\;y_1=b-d_2,\;\;y_2=b,\;\;y_3=b+d_2, \label{tep5subs1}
\end{equation}
where $a,\,b,\,d_1,\,d_2$, are arbitrary parameters, and with these values, \eqref{tep5eq2} and \eqref{tep5eq4} reduce to the two equations,
\begin{align}
3a^2+2d_1^2 &= 3b^2+2d_2^2, \label{tep5eq2a}\\
3a^4+12a^2d_1^2+2d_1^4 &= 3b^4+12b^2d_2^2+2d_2^4.\label{tep5eq4a}
\end{align}

Now \eqref{tep5eq2a} may be written as $3(a-b)(a+b)=-2(d_1-d_2)(d_1+d_2)$, and thus, its complete solution is readily obtained by writing,
\begin{equation}
a-b=2pr,\;\; 3(a+b)=12qs,\;\; d_1-d_2=2ps,\;\;-2(d_1+d_2)=12qr,
\label{tep5eq2b}
\end{equation}
and is given by,
\begin{equation}
a = pr+2qs, \;\; b = -pr+2qs,\;\;d_1 = ps-3qr, \;\; d_2 = -ps-3qr.
\label{tep5eq2c}
\end{equation}
where $p, \,q,\,r$ and $s$ are arbitrary parameters. 

With these values of $a,\,b,\,d_1$ and $d_2$, \eqref{tep5eq4a} reduces, on transposing all terms to one side and removing the factor $48pqrs$,  to
\begin{equation}
(2r^2-s^2)p^2-(9r^2-8s^2)q^2=0. \label{tep5eq4b}
\end{equation}
On writing $p=qs^2V/(2r^2-s^2),\;r=sU$, Eq.~\eqref{tep5eq4b} reduces to the quartic equation,
\begin{equation}
V^2=18U^4-25U^2+8.\label{tep5eq4c}
\end{equation}
Eq.~\eqref{tep5eq4c} is a quartic model of an elliptic curve that reduces, under the birational transformation,
\begin{equation}
\begin{aligned}
U &= (9X-Y-104)/(11X-Y-236), \\
V &= (X^3-198X^2+916X+980Y+34336)/(11X-Y-236)^2,\\
X &= 2(14U^2-17U+V+4)/(U-1)^2,\\
 Y &= 2(36U^3-25U^2+11UV-25U-9V+16)/(U-1)^3,
\end{aligned}
\label{tep5birat}
\end{equation}
to the Weierstrass form of an elliptic curve given by 
\begin{equation}
Y^2=X^3-X^2-784X+8704. \label{tep5eq4d}
\end{equation}

It is readily determined using APECS (a software package written in MAPLE for working with elliptic curves) that \eqref{tep5eq4d} is an elliptic curve of rank 1, its Mordell-Weil basis being given by the rational point $P$ with co-ordinates $(X,\,Y)=(-8,\,120)$. We can thus find infinitely many rational points on the curve \eqref{tep5eq4d}, and working backwards, we can find infinitely many nontrivial solutions of the simultaneous equations \eqref{tep5eq2} and \eqref{tep5eq4}, and hence also of the simultaneous equations \eqref{tep5s6} and \eqref{tep5condsym}, in which the integers $x_i, \,i=1,\,2,\,3$, and $y_i, \,i=1,\,2,\,3$, are   in arithmetic progression.

While the rational point $P$ on the elliptic curve \eqref{tep5eq4d} corresponds to a trivial solution of the equations \eqref{tep5eq2} and \eqref{tep5eq4}, the points $2P$ and $3P$ given by $(569/25,\, -5772/125)$  and  $(9121912/591361,\,$ $2979279240/454756609)$ respectively, yield the following two nontrivial solutions of Eqs.~\eqref{tep5eq2} and \eqref{tep5eq4}: 
\begin{equation}  
(x_1,\,x_2,\,x_3,\,y_1,\,y_2,\,y_3)=(1965,\, 1121,\,277,\, 1025,\, -477,\,-1979); \label{tep5numex1}
\end{equation}
\begin{multline}
(x_1,\,x_2,\,x_3,\,y_1,\,y_2,\,y_3)= (-201642299,\,47046243,\,295734785,\\
299528843,\,187147999,\,74767155) .\label{tep5numex2}
\end{multline}

\subsubsection{} 
 To obtain  parametric  ideal symmetric solutions of the TEP of degree 5, we will first find a solution of \eqref{tepsk} with $k=3$  taking the numbers $x_i$  as the terms of  the two arithmetic progressions $ [a_j,\,n_j,\,d_1],\; j=1,\,2$,
and  the numbers  $y_i$  as the terms of the arithmetic progressions  $ [b_j,\,n_j,\,d_2],\; j=1,\,2$.
 We thus have to solve the following equations:

\begin{align}
S_1(a_1,\,n_1,\,d_1)+S_1(a_2,\,n_2,\,d_1)&=S_1(b_1,\,n_1,\,d_2)+S_1(b_2,\,n_2,\,d_2), \label{2tep5symsys1}\\
S_2(a_1,\,n_1,\,d_1)+S_2(a_2,\,n_2,\,d_1)&=S_2(b_1,\,n_1,\,d_2)+S_2(b_2,\,n_2,\,d_2), \label{2tep5symsys2}\\
S_3(a_1,\,n_1,\,d_1)+S_3(a_2,\,n_2,\,d_1)&=S_3(b_1,\,n_1,\,d_2)+S_3(b_2,\,n_2,\,d_2), \label{2tep5symsys3}
\end{align}
In addition,  we   impose the condition that 
\begin{equation}
a_1+(2n_1-1)d_1 +2d_1=a_2-(2n_2-1)d_1, \label{2tep5cond1}
\end{equation}
so that  the numbers $x_i$ actually constitute just a single arithmetic progression with common difference $2d_1$. 

The symmetric diophantine equations Eqs.~\eqref{2tep5symsys1}, \eqref{2tep5symsys2} are readily solved together with the linear equation \eqref{2tep5cond1} and their  solution is given by,
\begin{equation}
\begin{aligned}
a_1 &= 8n_1^3+6t^2n_1^2n_2+6t^2n_1n_2^2+8n_2^3-2n_1-2n_2+a_2,\\
b_1 &= 8n_1^3+2t(3t-4)n_1^2n_2-2t(3t-4)n_1n_2^2\\
& \quad \quad -(8t-8)n_2^3-2n_1+(2t-2)n_2+a_2,\\
 b_2 &= 8tn_1^3+4(3t-2)tn_1^2n_2+8tn_1n_2^2-2tn_1+a_2, \\
d_1 &= -4n_1^2-(3t^2-4)n_1n_2-4n_2^2+1,\\
 d_2 &= 4n_1^2-(3t^2-12t+4)n_1n_2+4n_2^2-1,
\end{aligned} \label{2tep5solabd}
\end{equation}
where $t$ is an arbitrary parameter.

On substituting the values of $a_j,\,b_j,\,d_j$ given by \eqref{2tep5solabd} in \eqref{2tep5symsys3}, we get the following condition:
\begin{multline}
16n_1n_2(n_1+n_2)t(t-2)(2n_1+2n_2+1)(2n_1+2n_2-1)\\
\quad \times (4n_1^2+6tn_1n_2-4n_1n_2+4n_2^2-1)(n_1-n_2)\\
\quad \times \{(4t-4)n_1^2+(3t^2-4t+4)n_1n_2+(t-1)(4n_2^2-1)\}=0. \label{2tep5cond2}
\end{multline}

Equating to 0 any factor on the left-hand side of \eqref{2tep5cond2}, except the last two factors, leads to a trivial result. Equating to 0 the second  last factor of \eqref{2tep5cond2}, we get $n_2=n_1$, and we now have a solution of the diophantine system \eqref{tepsk} with $k=3$, and on    applying  Lemma~\ref{lemTarry} twice in succession, taking $h=2d_1$ and $h=2d_2$ respectively, we  obtain a symmetric solution of the diophantine system \eqref{tep5s6}. This symmetric solution simplifies further on writing $t=(2n_1-1)p/(3n_1q)$, and the reduced form of this solution may be written as, 
\begin{equation}
\begin{aligned}
x_1& =  (2p^2+6q^2)n_1^2-(p^2+6pq-3q^2)n_1-p^2-3q^2,\\
 x_2& = 2(p+q)(p-3q)n_1^2-(3p^2+9q^2)n_1+(p+q)(p-3q),\\
 x_3& = 8pqn_1^2+(p^2+3q^2)n_1-(p+3q)(p-q),\\
x_4&=-x_3, \quad x_5=-x_2, \quad x_6=-x_1, \\
 y_1 &= (2p^2+6q^2)n_1^2-3(p+q)(p-3q)n_1+p^2+3q^2,\\
 y_2& =  8pqn_1^2-(p^2+3q^2)n_1+(p+q)(p-3q),\\
 y_3 &= 2(p+q)(p-3q)n_1^2-(p^2+3q^2)n_1-(p+3q)(p-q),\\
y_4&=-y_3, \quad y_5=-y_2, \quad y_6=-y_1,
\end{aligned} \label{2tep5sol1}
\end{equation}
where $n_1,\,p$ and $q$ are arbitrary parameters.

Next, we equate to 0 the last factor of \eqref{2tep5cond2} and to solve this quadratic equation in $n_2$, we equate its discriminant to a perfect square. We thus obtain the following values of $n_1$ and $n_2$:
\begin{equation}
\begin{aligned}
n_1&= -8m(t-1)(m^2-9t^4+24t^3+24t^2-96t+48)^{-1},\\
n_2&=-\{m^2-(6t^2-8t+8)m+3(t+2)(3t-2)(t-2)^2\}\\
& \quad \quad \times (m^2-9t^4+24t^3+24t^2-96t+48)^{-1},
\end{aligned}
\label{2tep5valn}
\end{equation}
where $m$ is an arbitrary parameter. 

We now have  a second  solution of the diophantine system \eqref{tepsk} with $k=3$. We again apply Lemma~\ref{lemTarry} twice in succession, taking $h=2d_1$ and $h=2d_2$ respectively, and obtain another  symmetric solution of the diophantine system \eqref{tep5s6} which, expressed in reduced form, is given by, 
\begin{equation}
\begin{aligned}
x_1 &= 2m^2-6t(t-2)m+6(t-1)(t+2)(3t-2),\\
 x_2 &= m^2t+16(t-1)m-3t(t+2)(3t-2), \\
x_3 &= (2t-2)m^2-2(3t^2-4t+4)m-6(t+2)(3t-2),\\
 x_4&=-x_3, \quad x_5=-x_2, \quad x_6=-x_1, \\
y_1 &= (2t-2)m^2-6t(t-2)m+6(t+2)(3t-2), \\
 y_2 &= m^2t-16(t-1)m-3t(t+2)(3t-2), \\
 y_3 &= 2m^2+2(3t^2-4t+4)m-6(t-1)(t+2)(3t-2),\\
y_4&=-y_3, \quad y_5=-y_2, \quad y_6=-y_1,
\end{aligned} \label{2tep5sol2}
\end{equation}
where $m$ and $t$ are arbitrary parameters.

As a numerical example, taking $m=1,\,t=3$ in \eqref{2tep5sol2} yields the solution, 
\begin{equation}
\pm 101,\, \pm 70,\, \pm 61 \stackrel{5}{=} \pm 49,\, \pm 86,\, \pm 95.
\end{equation}

The solutions \eqref{2tep5sol1} and \eqref{2tep5sol2} of the diophantine system \eqref{tep5s6} immediately provide solutions of the simultaneous equations \eqref{tep5eq2} and \eqref{tep5eq4}. These solutions do not  satisfy the conditions \eqref{lincondneq}.

\subsection{Nonsymmetric solutions of the TEP of degree 5} 

\hspace{0.25in} Apart from a finite number of numerical solutions, only one parametric ideal nonsymmetric solution of the TEP of degree 5 in terms of polynomials of degree 11 has been published \cite{Cho2}.   We will now obtain such a parametric  ideal nonsymmetric solution  in terms of polynomials of degree 10 and show how infinitely many nonsymmetric solutions of \eqref{tep5s6}  may be obtained. 

While we can obtain  nonsymmetric solutions of \eqref{tep5s6} by our general method, such solutions can also be obtained by imposing the condition that  the solution \eqref{2tep4valXY} of the diophantine system \eqref{tep4s6} also satisfies the relation,
\begin{equation}
 \sum_{i=1}^6X_i^5=\sum_{i=1}^6Y_i^5. \label{condnstep5}
\end{equation}

Now \eqref{condnstep5} reduces to the following condition,
\begin{equation}
\begin{aligned}
4d^2&=(27f^4-14f^2g^2+3g^4)u^4-32fg(3f^2-g^2)u^3v\\
& \quad \;\;-(126f^4-108f^2g^2+14g^4)u^2v^2+96fg(3f^2-g^2)uv^3\\
& \quad \;\; +(243f^4-126f^2g^2+27g^4)v^4. 
\end{aligned} \label{condnstep5a}
\end{equation}

 The quartic function of $u$ and $v$ on the right-hand side of \eqref{condnstep5a}    becomes a perfect square 
 when $u=\pm v$ or $u=\pm 3v$. While these values of $u$ and $v$ lead to trivial results, we can obtain infinitely many values of $u,\,v$ that make the quartic function  on the right-hand side of \eqref{condnstep5a} a perfect square by following the  method described by Fermat \cite[p.\  639]{Di2}, one such solution being  $u=-3(3f^2-g^2),\;\; v=3f^2+8fg-g^2$. With these values of $u$ and $v$, Eq.~\eqref{condnstep5a} can be solved to get a rational value of $d$, and we thus obtain a nonsymmetric of \eqref{tep5s6} in terms of two arbitrary parameters $f$ and $g$. Denoting the polynomial $c_0f^n+c_1f^{n-1}g+c_2f^{n-2}g^2+\cdots+c_ng^n$ by $[c_0,\,c_1,\,c_2,\,\ldots,\,c_n]$, this two-parameter solution may be written in reduced form as follows:

\begin{equation}
\begin{aligned}
x_1&= [1701, 3888, 459, 1656, -5310, 2504, -1482, 616, 25, -24, -1],\\
x_2&= [243, 1944, 4509, -5256, -1314, -3112, 42, 40, -17, 48, -7],\\
x_3&= [-1215, -972, 2403, 1872, 4770, 4088, 798, 208, 157, -12, -1],\\
x_4&= [-1215, -2916, 459, 1656, 4122, -832, 1878, -248, -59, 36, -1],\\
x_5&= [243, 972, -3591, -936, 126, -1600, 354, -728, -17, -12, 5],\\
x_6&= [243, -2916, -4239, 1008, -2394, -1048, -1590, 112, -89, -36, 5],\\
y_1&= [243, -1944, -675, 5544, 4446, 2504, 1770, 184, -17, 48, -7],\\
y_2&= [-1215, -972, 459, -6552, -1062, -1600, -42, -104, 133, 12, -1],\\
y_3&= [243, -972, -4239, 1872, -2394, 4088, -1590, 208, -89, -12, 5],\\
y_4&= [243, 2916, 1593, -2232, -5634, -832, -1374, 184, -17, -36, 5],\\
y_5&= [1701, 3888, 459, 360, -126, -3112, 438, -584, -167, 24, -1],\\
y_6&= [-1215, -2916, 2403, 1008, 4770, -1048, 798, 112, 157, -36, -1].
\end{aligned} \label{solnstep6}
\end{equation}

As a numerical example, taking $f=2,\,g=-1$ in \eqref{solnstep6} yields the solution,
\begin{multline}
-87973,\,121805,\,-20525,\,52947,\,-108623,\,42369 \\
\stackrel{5}{=} 65869,\, 21507, \, -98863, \, -100895,\, -8325,\, 120707. \label{solnstep6ex1}
\end{multline}

As we can obtain infinitely many values of $u,\,v$ such that the quartic function  on the right-hand side of \eqref{condnstep5a}  becomes a perfect square, we can obtain infinitely many parametric nonsymmetric solutions of the TEP of degree 5.

\section{Ideal solutions of the Tarry-Escott problem of degree 6}\label{tepdeg6}
  \setcounter{equation}{0}

\hspace{0.25in} In this Section we will obtain  a parametric ideal solution of the TEP of degree 6, that is, of the diophantine system,
\begin{equation}
\sum_{i=1}^7x_i^r=\sum_{i=1}^7y_i^r,\;\;\;r=1,\,2,\,\dots,\,6.
\label{tep6s7}
\end{equation}
and we will show how more parametric solutions may be obtained. We note that multi-parameter solutions of \eqref{tep6s7} have been given by Choudhry \cite[p.\  305]{Cho0} and Gloden \cite[p.\  43]{Glo}.

We  first find a solution of \eqref{tepsk} with  $k=5$ in which the numbers $x_i$ on the left-hand side of \eqref{tepsk} are the terms of the four arithmetic progressions $[a_1,\,n_1,\,d], [-a_1,\,n_1,\,d]$, $[a_2,\,n_2,\,d]$ and $[-a_2,\,n_2,\,d]$, while the numbers  $y_i$ on the right-hand side of \eqref{tepsk} are the terms of the four arithmetic progressions  $[b_1,\,n_1,\,d], [-b_1,\,n_1,\,d], [b_2,\,n_2,\,d]$ and $[-b_2,\,n_2,\,d]$.

With the above choice of $x_i,\,y_i$, it is clear that \eqref{tepsk} is identically true for  $r=1,\,3$ and 5. We thus have to solve just the following two equations obtained by taking $r=2$ and $r=4$ in \eqref{tepsk} respectively:
\begin{equation}
2n_1a_1^2+2n_2a_2^2 = 2n_1b_1^2+2n_2b_2^2, \label{2tep6eq2}
\end{equation}
\begin{multline}
2n_1a_1^4+2n_2a_2^4+\{4n_1(4n_1^2-1)a_1^2+4n_2(4n_2^2-1)a_2^2\}d^2\\
 =2n_1b_1^4+2n_2b_2^4+\{4n_1(4n_1^2-1)b_1^2+4n_2(4n_2^2-1)b_2^2\}d^2. \label{2tep6eq4}
\end{multline}

Now \eqref{2tep6eq2} may be written as,
\begin{equation}
n_1(a_1-b_1)(a_1+b_1)=-n_2(a_2-b_2)(a_2+b_2),
\end{equation}
and its complete solution obtained by writing,
\begin{equation}
\begin{aligned}
n_1(a_1-b_1)&=2n_1pr, \quad &a_1+b_1&=2n_2qs,\\
-n_2(a_2-b_2)&=2n_2ps,&a_2+b_2&=2n_1qr,
\end{aligned}
\end{equation}
is given in terms of arbitrary parameters $p,\,q,\,r,\, s$ by,
\begin{equation}
 a_1 = pr+n_2qs,\;\; a_2 = -ps+n_1qr,\;\; b_1 = -pr+n_2qs,\;\; b_2 = ps+n_1qr. \label{sol2tep6eq2}
\end{equation}

With  these values of $a_i,\,b_i,\;i=1,\,2$, Eq.~\eqref{2tep6eq4} reduces to the following equation of degree two in $r,\,s$ and $d$:
\begin{equation}
(p^2-n_1^2q^2)r^2-(p^2-n_2^2q^2)s^2+4(n_1^2-n_2^2)d^2=0. \label{2tep6eq4a}
\end{equation}

Now \eqref{2tep6eq4a} has a solution $r=2,\,s=2,\,d=q,$ and so its complete solution is easily obtained, and is given by
\begin{equation}
 \begin{aligned}
r&=2(p^2-n_1^2q^2)u^2-4(p^2-n_2^2q^2)uv+2(p^2-n_2^2q^2)v^2,\\
s&=-2(p^2-n_1^2q^2)u^2-(p^2-n_1^2q^2)uv-2(p^2-n_2^2q^2)v^2,\\
d&=-(p^2-n_1^2q^2)qu^2+(p^2-n_2^2q^2)qv^2,\\
\end{aligned}
\label{2tep6eq4b}
\end{equation}
where $u,\,v$ are arbitrary parameters. 

Substituting the values of $r$ and $s$  given by \eqref{2tep6eq4b} in \eqref{sol2tep6eq2}, we get the values of $a_1,\,a_2,\,b_1$ and $b_2$  in terms of the parameters $p,\,q,\,n_1,\,n_2,\,u$ and $v$. We now have a solution of \eqref{tepsk} with $k=5$, and on applying Lemma~\ref{lemTarry} taking $h=2d$, we get a  solution of the diophantine system \eqref{tepsk} with $k=6$ and $s=8$ 
in terms of the parameters $p,\,q,\,n_1,\,n_2,\,u$ and $v$. A pair of terms, one on each side of this solution,  cancels out if the following condition is satisfied:
\begin{multline}
(p^2-n_1^2q^2)\{p+(n_1-n_2)q\}u^2-(2p^3-2n_2p^2q-2n_2^2pq^2\\
+2n_1^2n_2q^3)uv+(p^2-n_2^2q^2)\{p-(n_1+n_2)q\}v^2=0. \label{canc17}
\end{multline}

Taking $p=-(n_1-n_2)q$, the coefficient of $u^2$ in \eqref{canc17} vanishes, and we get $u=n_1^2(n_1-2n_2),\;v=n_1^3-3n_1^2n_2+n_2^3$, as a solution of \eqref{canc17}, and finally we obtain, after cancelling out a pair of terms,  a symmetric solution of \eqref{tep6s7} which may be written in the reduced form as follows:
\begin{equation}
\begin{aligned}
x_1&=-4n_1(n_1-n_2)(n_1+n_2)(n_1^2-3n_1n_2+n_2^2),\\
x_2&= -4n_1(n_1^4-4n_1^3n_2+5n_1^2n_2^2-n_2^4),\\
x_3&= 4(n_1^2-n_1n_2+n_2^2)(n_1^3-3n_1^2n_2+n_2^3),\\
x_4&= -4n_2(n_1^4-4n_1^3n_2+n_1^2n_2^2+2n_1n_2^3-n_2^4),\\
x_5&= -4n_1n_2(n_1-2n_2)(n_1^2+n_1n_2-n_2^2), \\
x_6&=4(n_1-n_2)(n_1^4-2n_1^3n_2-n_1^2n_2^2+n_2^4),\\
x_7&= 4n_2(2n_1-n_2)(n_1-n_2)(n_1^2-n_1n_2-n_2^2),\\
y_i&=-x_i,\;i=1,\,2,\,\ldots,\,7.
\end{aligned} \label{sol2tep6}
\end{equation} 
where $n_1$ and $n_2$ are arbitrary parameters. 

We note that the condition \eqref{canc17} is a quadratic  equation in $u,\,v$, and its discriminant is a quartic function of $p$ and $q$. We have already found one pair of values of $p$ and $q$ that make the discriminant a perfect square. Thus, following the method described by Fermat \cite[p.\  639]{Di2}, we can find infinitely many values of $p,\,q$ such that the discriminant becomes a perfect square, and hence we can obtain infinitely many solutions of \eqref{canc17}. We can thus find infinitely many  parametric ideal solutions of the TEP of degree 6.  

As a numerical example, taking $n_1=3,\,n_2=1$ in \eqref{sol2tep6} yields the solution,
\[
-66,\, -134,\, 133,\, 47,\, 8,\, 87,\, -75 \stackrel{6}{=} 66,\, 134,\, -133,\, -47,\, -8,\, -87,\, 75.
\]

\section{Ideal solutions of the Tarry-Escott problem of degree 7}\label{tepdeg7}
 \setcounter{equation}{0}
\hspace{0.25in}  We will now obtain ideal solutions of the  TEP of degree 7, that is, of the diophantine system,
 \begin{equation}
\sum_{i=1}^8x_i^r=\sum_{i=1}^8y_i^r,\;\;\;r=1,\,2,\,\dots,\,7.
\label{tep7s8}
\end{equation}
Till now only one parametric solution of \eqref{tep7s8} has been published \cite{Che}. This is a symmetric solution that  satisfies the additional conditions,
\begin{equation}
\begin{aligned}
x_5&=-x_4,\; &x_6&=-x_3,\;\; &x_7&=-x_2,\;\; &x_8&=-x_1,\\
y_5&=-y_4,\; &y_6&=-y_3,\;\; &y_7&=-y_2,\;\; &y_8&=-y_1. \label{tep7condsym}
\end{aligned}
\end{equation}

We obtain infinitely many numerical solutions, as well as a parametric solution of \eqref{tep7s8}, and show how infinitely many parametric solutions may be obtained. All the solutions that we obtain are symmetric, and hence they also provide solutions of the diophantine system,
\begin{equation}
\sum_{i=1}^4x_i^r=\sum_{i=1}^4y_i^r,\;\;\;r=2,\,4,\,6.
\label{tep246}
\end{equation}

\subsection{}
In Section~4.1.1 we have  described a method of obtaining infinitely many numerical solutions  of the simultaneous equations \eqref{tep5s6} and \eqref{tep5condsym} in which the integers 
$x_i,\,i=1,\,2,\,3$, and the integers $y_i,\,i=1,\,2,\,3$, are the terms of two arithmetic progressions with common differences $d_1$ and $d_2$ such that $d_1 \neq d_2$. Applying Lemma~\ref{lemTarry}  twice, in succession, to  such solutions taking $h=d_1$ and $h=d_2$ respectively,  immediately yields, on cancellation of common terms on either side, infinitely many  symmetric solutions of \eqref{tep7s8}. As an example,  the two numerical solutions  \eqref{tep5numex1} and \eqref{tep5numex2} yield 
the following two solutions of \eqref{tep7s8}:
\begin{equation}
\pm 448,\, \pm 677,\, \pm 1154,\, \pm 1569 \stackrel{7}{=} \pm 303,\, \pm 818,\, \pm 1099,\, \pm 1576;
\end{equation}

\begin{multline}
\pm 181944317,\,  \pm 134898074,\,  \pm 240031768,\,  \pm 52883769 \\
 \stackrel{7}{=}  \pm 238134739,\, \pm 191088496,\, \pm 115687497,\, \pm 71460502.
\end{multline}

\subsection{}
We will now obtain a parametric solution of the diophantine system \eqref{tep7s8}. We  first find a solution of \eqref{tepsk} with $k=5$ taking the numbers $x_i$ on the left-hand side as the terms of the  two arithmetic progressions $[a,\,n,\,d_1]$ and $[-a,\,n,\,d_1]$ and the  numbers $y_i$ as the terms of the two  arithmetic progressions $[b,\,n,\,d_2]$ and $[-b,\,n,\,d_2]$.

With the above  values of   $x_i,\,y_i$, it is clear that \eqref{tepsk} is identically true for $r=1,\,3$ and 5. We thus have to solve just the following two equations obtained by taking $r=2$ and $r=4$ in \eqref{tepsk} respectively:
\begin{equation}
2na^2+(2/3)n(4n^2-1)d_1^2 = 2nb^2+(2/3)n(4n^2-1)d_2^2, \label{tep7eq2}\\
\end{equation}
and
\begin{multline}
2na^4+4n(4n^2-1)a^2d_1^2+(2/15)n(4n^2-1)(12n^2-7)d_1^4\\
 = 2nb^4+4n(4n^2-1)b^2d_2^2+(2/15)n(4n^2-1)(12n^2-7)d_2^4. \label{tep7eq4}
\end{multline}

 Now \eqref{tep7eq2} may be written as,
\begin{equation}
2n(a-b)(a+b) = -(2/3)n(2n-1)(2n+1)(d_1-d_2)(d_1+d_2),\label{tep7eq2a}
\end{equation}
and its complete solution obtained by writing,
\begin{equation}
\begin{aligned}
a-b&=2(2n-1)pr,&\quad a+b&=2(2n+1)qs,\\
-(d_1-d_2)/3&=2ps,&\quad d_1+d_2&=2qr,
\end{aligned}
\end{equation}
is given  in terms of arbitrary parameters $p,\,q,\,r,\, s$ by,
\begin{equation}
\begin{aligned}
a &= (2n-1)pr+(2n+1)qs,\quad & b& = -(2n-1)pr+(2n+1)qs, \\
d_1& = -3ps+qr,& d_2& = 3ps+qr.
\end{aligned}
 \label{soltep7eq2}
\end{equation}

Using  these values of $a,\,b,\,d_1$, and $d_2$, \eqref{tep7eq4} reduces, on removing the factor $pqrsn(2n-1)(2n+1)$,  to the following equation: 
\begin{equation}
\{5(2n-1)^2r^2-9(4n^2+1)s^2\}p^2-\{(4n^2+1)r^2-5(2n+1)^2s^2\}q^2=0. \label{tep7eq4a}
\end{equation}

Now \eqref{tep7eq4a} will have a rational solution for $p$ and $q$ if and only if  the following quartic function in $r$ and $s$,
\begin{equation}
\{5(2n-1)^2r^2-9(4n^2+1)s^2\}\{(4n^2+1)r^2-5(2n+1)^2s^2\}, \label{tep7dis}
\end{equation}
 becomes a perfect square. We observe that the quartic function \eqref{tep7dis}  does become a perfect square when $r=s$. While this leads to a trivial solution, we can use this solution to find  values of $r$ and $s$ that make the quartic function \eqref{tep7dis}  a perfect square by the  method described by Fermat \cite[p.\ 639]{Di2}. We thus get,
\begin{equation}
r=-(4n^2+9n+1),\quad s=4n^2+n+1, \label{tep7rs}
\end{equation}
and with these values of $r$ and $s$, Eq.~\eqref{tep7eq4a} has the following solution:
\begin{equation}
p = 8n^3-6n^2+3n-1,\quad q = -(8n^3-14n^2-7n+1). \label{tep7pq}
\end{equation}

Substituting the values of $p,\,q,\,r$ and $s$ given by \eqref{tep7rs} and \eqref{tep7pq} in Eq.~\eqref{soltep7eq2}, we obtain the following solution of equations \eqref{tep7eq2} and \eqref{tep7eq4}:
\begin{equation}
\begin{aligned}
a &= -2n(16n^4-17n^2-3),\quad & b& = 48n^4+17n^2-1, \\
d_1& = -16n^4-47n^2-1, &d_2& = 32n^4-26n^2+2. 
\end{aligned}
\label{soltep7eq24}
\end{equation}

We now have a solution of \eqref{tepsk} with $k=5$ and  on applying Lemma~\ref{lemTarry} twice in succession, taking $h=2d_1$ and $h=2d_2$ respectively, we get a symmetric solution of  the diophantine system \eqref{tep7s8}. The reduced form of this solution may be written as follows: 
\begin{equation}
\begin{aligned}
x_1&=16n^4-64n^3-13n^2-4n+1, \\
x_2 &=32n^5-16n^4+30n^3+13n^2-2n-1, \\
x_3 &= -32n^5+16n^4+26n^3-15n^2-2n-1,\\
 x_4 &=-32n^5-32n^4+26n^3-32n^2-2n, \\
x_5&=-x_4,\quad x_6=-x_3,\quad x_7=-x_2,\quad x_8=-x_1,\\
y_1 &=-32n^5+32n^4+26n^3+32n^2-2n,\\
y_2 &= -32n^5-16n^4+26n^3+15n^2-2n+1, \\
y_3 &=16n^4+64n^3-13n^2+4n+1, \\
y_4 &= 32n^5+16n^4+30n^3-13n^2-2n+1,\\ 
y_5&=-y_4,\quad y_6=-y_3,\quad y_7=-y_2,\quad y_8=-y_1,
\end{aligned} \label{soltep7}
\end{equation}
where $n$ is an arbitrary parameter.

We note that the quartic function \eqref{tep7dis} can be made a perfect square for infinitely many values of $r$ and $s$ that may be obtained by repeated application of  Fermat's method mentioned above. We can thus obtain infinitely many parametric ideal  solutions of  the TEP of degree 7.

As a numerical example, taking $n=2$ in \eqref{soltep7} yields the solution,
\begin{equation}
\pm 63,\, \pm 211, \pm 125, \pm 292 \stackrel{7}{=} \pm 36,\, \pm 203, \pm 145, \pm 293.
\end{equation}

\section{Some diophantine systems related to the Tarry-Escott problem}\label{relsys}
\hspace{0.25in} In this Section we briefly consider some diophantine systems that are closely related to the TEP. We can  obtain new solutions of several such diophantine systems using the new approach to   the TEP described in this paper. We restrict ourselves to giving a few examples.
\setcounter{equation}{0}

\subsection{}
In this subsection, we  consider the following diophantine system,
\begin{equation}
\sum_{i=1}^{k+1}x_i^r=\sum_{i=1}^{k+1}y_i^r,\;\;\;r=1,\,2,\,\dots,\,k,\,k+2.
\label{tepskaug}
\end{equation}
To obtain solutions of this diophantine system, we will use the following theorem proved by Gloden \cite[p.\  24]{Glo}.

\begin{thm}\label{thmGloden}
 If there exist integers $x_i,\,y_i,\,i=1,\,2,\,\ldots,\,k+1,$ such that the relations \eqref{tepsk} are satisfied with $s=k+1$, then
\begin{equation}
\sum_{i=1}^{k+1}(x_i+d)^r=\sum_{i=1}^{k+1}(y_i+d)^r,\;\;\;r=1,\,2,\,\dots,\,k,\,k+2,
\label{teplemGlo}
\end{equation}
where
\begin{equation}
d=-\left(\sum_{i=1}^{k+1}x_i\right)/(k+1).
\end{equation}
\end{thm}

It is to be noted that  Theorem~\ref{thmGloden} yields nontrivial solutions of the diophantine system \eqref{tepskaug} only when we apply it  to  ideal nonsymmetric solutions.   

 We have already explicitly obtained two parametric ideal nonsymmetric solutions of the TEP of degree 4,  the first solution being  given by \eqref{sol3tep4s5x}and  \eqref{sol3tep4s5y} and the second one by \eqref{tep4valxy}. We further note that these solutions are already in reduced form and  thus each of them satisfies the additional condition $\sum_{i=1}^5x_i=0$. It now immediately follows from Theorem~\ref{thmGloden} that  these two parametric solutions  as well as the numerical solutions \eqref{tep4ex1} and \eqref{tep4ex2}  also satisfy the diophantine system,
\begin{equation}
\sum_{i=1}^5x_i^r=\sum_{i=1}^5y_i^r,\;\;\;r=1,\,2,\,3,\,4,\,6.
\end{equation}

Similarly, it follows from Theorem~\ref{thmGloden} that the parametric ideal nonsymmetric solution of the TEP of degree 5 given by \eqref{solnstep6} and the numerical solution \eqref{solnstep6ex1}  are also  solutions of the diophantine system,
\begin{equation}
\sum_{i=1}^6x_i^r=\sum_{i=1}^6y_i^r,\;\;\;r=1,\,2,\,3,\,4,\,5,\,7.
\end{equation}

\subsection{}
 In this subsection,  we  consider the following diophantine system,
\begin{equation}
\begin{aligned}
\sum_{i=1}^sx_i^r&=\sum_{i=1}^sy_i^r,\;\;\;r=1,\,2,\,\dots,\,k,\\
\prod_{i=1}^sx_i&=\prod_{i=1}^sy_i.
\end{aligned} \label{tepskeqprod}
\end{equation}
A detailed discussion of this diophantine system is given in \cite{Cho4} where it is  shown  that for the existence of a nontrivial solution of the diophantine system \eqref{tepskeqprod}, it is necessary that $s \geq k+2.$ Here we restrict ourselves to finding solutions of the diophantine system \eqref{tepskeqprod} when $(k,\,s)=(4,\,6)$ and also when $(k,\,s)=(5,\,7)$ by applying a lemma proved by Choudhry \cite[Lemma 3, pp.\  766-767]{Cho4}.

\subsubsection{}

Applying the aforesaid lemma  to either of the two solutions \eqref{sol3tep4XY} and \eqref{2tep4valXY}   of the diophantine system \eqref{tep4s6}, we immediately get two multi-parameter  solutions of    \eqref{tepskeqprod} with $(k,\,s)=(4,\,6)$. As these solutions are cumbersome to write, we omit writing them explicitly and give 
below just one three-parameter solution derived from the solution \eqref{2tep4valXY} in which we have taken $u=1,\,v=1$. This solution is as follows:
\begin{equation}
\begin{aligned}
x_1&=(6f^2+12fg+6g^2+d)(-30f^2+4fg+2g^2+d), \\
x_2&=(-6f^2+4fg+10g^2+d)(30f^2-4fg-2g^2+d),\\
  x_3&=(18f^2+20fg+2g^2+d)(-18f^2+12fg-2g^2+d), \\
x_4&=	(6f^2-4fg-10g^2+d)(18f^2-12fg+2g^2+d),\\
x_5&= (-6f^2+20fg-6g^2+d)(-18f^2-20fg-2g^2+d),\\
x_6&= (6f^2-20fg+6g^2+d)(-6f^2-12fg-6g^2+d),\\
y_1&=(-18f^2+12fg-2g^2+d)(-6f^2+4fg+10g^2+d),\\
y_2&= (6f^2-20fg+6g^2+d)(18f^2+20fg+2g^2+d), \\
y_3&=(-6f^2+20fg-6g^2+d)(6f^2+12fg+6g^2+d),\\
y_4&= (30f^2-4fg-2g^2+d)(-6f^2-12fg-6g^2+d),\\
y_5&= (-30f^2+4fg+2g^2+d)(6f^2-4fg-10g^2+d),\\
y_6&= (18f^2-12fg+2g^2+d)(-18f^2-20fg-2g^2+d),
\end{aligned} \label{soltep4eqprod}
\end{equation}
where $f,\,g$ and $d$ are arbitrary parameters. While this solution is of degree 2 in the parameter $d$,  a solution of degree 4 has been published  earlier \cite{Cho4}.

As a numerical example, taking $f=2,\,g=1,\,d=1$ in \eqref{soltep4eqprod}, we get the solution,
\begin{align*}
5995,\, 555,\, 5635,\, -357,\, 1243,\, -477\stackrel{4}{=}-245,\, 1035,\, -605, \,5883,\, 763,\, 5763,\\
5995.555.5635.(-357).1243.(-477)=(-245).1035.(-605).5883.763.5763 .
\end{align*}

\subsubsection{}
We will now solve the diophantine system \eqref{tepskeqprod}  when $(k,\,s)=(5,\,7).$  

We will first find a solution of the diophantine system \eqref{tepsk} with $k=2$ and $s=2m_1+2m_2=2n$, taking the numbers $x_i$ as the terms of two arithmetic progressions $[a_j,\,m_j,\,d_j],\;j=1,\,2,$ and the the numbers $y_i$ as the terms of the arithmetic progression $[0,\,n,\,d_3]$, where $d_1,\,d_2$ and $d_3$ are distinct rational numbers. We now have to solve the following two equations obtained by taking $r=1$ and $r=2$ respectively in \eqref{tepsk}:
\begin{equation}
m_1a_1+n_2m_2=0, \label{prod2a}
\end{equation} 
\begin{multline}
 m_1a_1^2+m_1(4m_1^2-1)d_1^2/3+m_2a_2^2+m_2(4m_2^2-1)d_2^2/3\\
=n(4n^2-1)d_3^2/3. \label{prod2b}
\end{multline}
In addition, we impose the following auxiliary conditions so that when we apply Lemma~\ref{lemTarry} with $h=2d_3$ to the  solution of the diophantine system \eqref{tepsk} that is being obtained, the resulting solution of \eqref{tepsk} with $k=3$ will consist of the terms of just four arithmetic progressions.
\begin{align}
(2n+1)d_3 -\{a_1+(2m_1-1)d_1\}&=2d_1,\label{prod2c}\\
a_2-(2m_2-1)d_2+2d_3+ (2n-1)d_3&=2d_2. \label{prod2d}
\end{align}

We now solve the four equations \eqref{prod2a}, \eqref{prod2b}, \eqref{prod2c}, \eqref{prod2d} for $a_1,\,a_2,\,d_1,$ $d_2$ and $d_3$, and to the resulting solution of \eqref{tepsk} with $k=2$, we apply Lemma~\ref{lemTarry} three times in succession, taking $h=2d_3,\;h=2d_1$ and $h=2d_2$ respectively,   and thus obtain, after cancellation of common terms,  a solution of \eqref{tepsk} with $k=5$ and $s=8$. On taking $n_1=-3/4$, two more terms - one from each side - cancel out, and we obtain a solution of \eqref{tepsk} with $k=5$ and $s=7$. 

We now apply the aforesaid lemma proved by Choudhry \cite[Lemma 3, pp.\  766-767]{Cho4} to obtain a solution of the diophantine system \eqref{tepskeqprod} with $k=5$ and $s=7$ in terms of the parameter $m_2$. On replacing $m_2$ by $m$, this solution may be written as follows:
\begin{equation}
\begin{aligned}
x_1&=(4m+3)(8m-3),&x_2&=-2(2m+3)(12m-1),\\
x_3&=-4(3m+1)(4m-9),&x_4&=6(4m+1)(8m+1),\\
x_5&=3(4m+1)(16m-3),&x_6&=8(4m+3)(m-1),\\
x_7&=-4(16m+9)(2m-1),&y_1&=8m+1)(4m-9),\\
y_2&=(16m+9)(12m-1),&y_3&=4(2m+3)(4m+3),\\
y_4&=8(3m+1)(8m-3),&y_5&=-6(4m+1)(2m-1),\\
y_6&=-2(4m+1)(16m-3),&y_7&=-12(4m+3)(m-1).
\end{aligned}
\end{equation}

As a numerical example, taking $m=-1$ yields the solution,
\[
\begin{aligned}
11,\, 26,\, -104,\, 126,\, 171,\, 16,\, -84 &\stackrel{5}{=}91,\, 91,\, -4,\, 176,\, -54,\, -114,\, -24,\\
11.26.(-104).126.171.16.(-84) &=91.91.(-4).176.(-54).(-114).(-24).
\end{aligned}
\]

\section{Concluding Remarks} 
\hspace{0.25in} In this paper we have described a new method to derive solutions of the TEP, and obtained several new parametric ideal solutions of the TEP of degree $\leq 7$. The method can be applied to obtain many other new solutions of the TEP and related diophantine systems. It may be possible to apply the method described in this paper to obtain the  complete ideal solution of the  TEP of degrees 4 and 5  but we leave this as an open problem. Similarly it is perhaps possible to apply this method to obtain parametric ideal solutions of the TEP of degrees $\geq 8$ but this is also left for future investigations.

\section{Acknowledgements} {I am grateful to the Harish-Chandra Research Institute, Allahabad, India for providing me   all necessary facilities that have helped me  to pursue my research work in mathematics. I also wish to thank the referee for his comments and suggestions that have led to improvements in the paper.


\begin{thebibliography}{99}

\bibitem {Bor}   P. Borwein,   and C.  Ingalls, The Prouhet-Tarry-Escott
problem revisited,  {\it Enseign. Math.} {\bf 40}  (1994),
3--27.


\bibitem {Che}J. Chernick,  Ideal Solutions of the Tarry-Escott problem, {\it Amer.
Math. Monthly} {\bf 44} (1937), 626--633.


\bibitem{Cho0} A. Choudhry,  Symmetric diophantine systems, {\it Acta Arith.} {\bf 59}
(1991),  291--307.

\bibitem {Cho1}   A. Choudhry,  Ideal Solutions of the Tarry-Escott problem of
degree four and a related diophantine system, {\it Enseign. Math.} {\bf 46} (2000), 313--323.


\bibitem{Cho2}  A. Choudhry,  Ideal solutions of the Tarry-Escott problem of degrees four
and five and related diophantine systems, {\it Enseign. Math.}  {\bf 49} (2003),   101--108.

\bibitem{Cho3}  A. Choudhry, Ideal solutions of the Tarry-Escott problem of degree eleven
with applications to sums of thirteenth powers, {\it Hardy-Ramanujan J.} {\bf 31} (2008),
1--13.

\bibitem{Cho4}  A. Choudhry,  Equal Sums of Like Powers and Equal Products of Integers,
{\it Rocky Mountain J.  Math.}  {\bf 43} (2013),  763--792.


\bibitem{Di1} L. E. Dickson,   {\it Introduction to the theory of numbers,} Dover
Publications, New York, 1957, reprint.

\bibitem{Di2} L. E. Dickson,  {\it History of the theory of numbers},   Vol. 2,
Chelsea Publishing Company, New York, 1992, reprint.

\bibitem {Dor1}  H. L. Dorwart,   Sequences of ideal solutions in the Tarry-Escott
problem,  {\it Bull. Amer. Math. Soc.} {\bf 53} (1947), 381--391.

\bibitem {Dor2}  H. L. Dorwart,   and O. E. Brown, The Tarry-Escott problem,  {\it Amer. Math.
Monthly} {\bf  44} (1937), 613--626.

\bibitem  {Glo}  A. Gloden,   {\it Mehrgradige Gleichungen}, Noordhoff, Groningen, 1944.

\bibitem {Let}  A. Letac,  {\it Gazeta Mathematica} {\bf 48} (1942), 68--69.

\bibitem {Shu}  C. Shuwen,   Equal sums of like powers,  Website: http://euler.free.fr/eslp 

\bibitem {Smy}  C. J. Smyth, Ideal 9$^{th}$-order multigrades and Letac's elliptic
curve,  {\it Math. Comp.} {\bf 57} (1991), 817--823.
\end{thebibliography}
\end{document}